\documentclass[12pt]{article}

\usepackage{amsmath, epsfig, cite}
\usepackage{amssymb}
\usepackage{amsfonts}
\usepackage{latexsym}
\usepackage{float}

\newtheorem{thm}{Theorem}[section]

\newtheorem{conj}[thm]{Conjecture}
\newtheorem{lem}[thm]{Lemma}

\numberwithin{equation}{section}

\newcommand{\qed}{{\hfill$\square$}\medskip}

\begin{document}

\begin{center}
{\Large\bf Congruences on sums of super Catalan numbers
}
\end{center}

\vskip 2mm \centerline{Ji-Cai Liu}
\begin{center}
{\footnotesize Department of Mathematics, Wenzhou University, Wenzhou 325035, PR China\\
{\tt jcliu2016@gmail.com} }
\end{center}


\vskip 0.7cm \noindent{\bf Abstract.} In this paper, we prove two congruences on
the double sums of the super Catalan numbers (named by Gessel), which were recently conjectured by
Apagodu.


\vskip 2mm
\noindent{\it MR Subject Classifications}: 11A07, 05A19, 05A10
\section{Introduction}
It is well-known that the Catalan numbers
\begin{align*}
C_n=\frac{1}{n+1}{2n\choose n}
\end{align*}
are integers and occur in various counting problems. We refer to \cite{stanley-b-1999} for many different combinatorial interpretations of the Catalan numbers. The closely related central binomial coefficients are given by ${2n\choose n}$ for $n\in\mathbb{N}$.

Both Catalan numbers and central binomial coefficients possess many interesting arithmetic properties.
Sun and Tauraso \cite{st-ijnt-2011} proved that for primes $p\ge 5$,
\begin{align*}
&\sum_{k=0}^{p-1}{2k\choose k}\equiv \left(\frac{p}{3}\right)\pmod{p^2},\\
&\sum_{k=0}^{p-1}C_k\equiv\frac{3}{2}\left(\frac{p}{3}\right)- \frac{1}{2}\pmod{p^2},
\end{align*}
where $\left(\frac{\cdot}{p}\right)$ denotes the Legendre symbol. Recently, Mattarei and Tauraso \cite{mt-jnt-2018} showed that
\begin{align}
&\sum_{k=0}^{q-1}{2k\choose k}x^k\equiv (1-4x)^{\frac{q-1}{2}}\pmod{p},\label{mt-1}\\
&\sum_{k=0}^{q-1}C_kx^{k+1}\equiv\frac{1-(1-4x)^{\frac{q+1}{2}}}{2}-x^q\pmod{p},\label{mt-2}
\end{align}
where $q$ is a power of an odd prime $p$. For more congruence properties on these numbers we refer to \cite{gz-aam-2010, st-aam-2010, tauraso-aam-2012}.

In 1874, E. Catalan observed that the numbers
\begin{align*}
S(m,n)=\frac{{2m\choose m}{2n\choose n}}{{m+n\choose m}}
\end{align*}
are integers. Since $S(1,n)/2$ coincides with $C_n$, these numbers $S(m,n)$ are named {\it super Catalan numbers} by Gessel \cite{gessel-jsc-1992}. These numbers should not confused with the Schr\"oder--Hipparchus numbers, which are sometimes also called super Catalan numbers. Some interpretations of $S(m,n)$ for some special values of $m$ have been studied by several authors (see, e.g., \cite{ag-jis-2014, cw-a-2012, ps-rsa-2006}).
It is still an open problem to find a general combinatorial interpretation for the super Catalan numbers.

Our interest concerns the following two conjectures by Apagodu \cite[Conjecture 2]{apagodu-injt-2018}.
\begin{conj} (Apagodu)
For any odd prime $p$, we have
\begin{align}
\sum_{i=0}^{p-1}\sum_{j=0}^{p-1}S(i,j)\equiv \left(\frac{p}{3}\right)\pmod{p},\label{conj-1}\\[5pt]
\sum_{i=0}^{p-1}\sum_{j=0}^{p-1}(3i+3j+1)S(i,j)\equiv -7\left(\frac{p}{3}\right)\pmod{p}.\label{conj-2}
\end{align}
\end{conj}

In Section 2, we provide a proof of \eqref{conj-1} which makes use of a combinatorial identity.
\begin{thm} \label{t-1}
The congruence \eqref{conj-1} is true.
\end{thm}

We prove \eqref{conj-2} by establishing the following congruence.
\begin{thm} \label{t-2}
For any prime $p\ge 5$, we have
\begin{align}
\sum_{i=0}^{p-1}\sum_{j=0}^{p-1}(i+j)S(i,j)\equiv -\frac{8}{3}\left(\frac{p}{3}\right)\pmod{p}. \label{a-1}
\end{align}
\end{thm}
From \eqref{conj-1} and \eqref{a-1}, we deduce \eqref{conj-2} for $p\ge 5$.
It is routine to check that \eqref{conj-2} also holds for $p=3$.

\section{Proof of Theorem \ref{t-1}}
In order to prove Theorem \ref{t-1}, we need the following identity.
\begin{lem}
For any non-negative integer $n$, we have
\begin{align}
\sum_{i=0}^n\sum_{j=0}^n(-4)^{i+j}\frac{ {n\choose i}{n\choose j}}{ {i+j\choose i}}=
\frac{(-3)^n(2n-1)}{4(n+1)}+\frac{4^n}{{2n\choose n}}
\left(\frac{1}{2}-\sum_{k=0}^nC_k\left(-\frac{3}{4}\right)^{k+1}\right),\label{b-1}
\end{align}
where $C_k$ denotes the $k$th Catalan number.
\end{lem}
{\noindent \it Proof.}
Applying the multi-Zeilberger algorithm \cite{az-aam-2006}, we find that the left-hand side of \eqref{b-1} satisfies the recurrence:
\begin{align*}
-18(n+1)s(n)+3(2n-5)s(n+1)+2(5n+6)s(n+2)+(5+2n)s(n+3)=0.
\end{align*}
It is routine to check that the right-hand side of \eqref{b-1} also satisfies this recurrence and both sides of \eqref{b-1} are equal for $n=0,1,2$.
\qed

{\noindent \it Proof of \eqref{conj-1}.}
Let $n=\frac{p-1}{2}$. We split the double sum on the left-hand side of \eqref{conj-1} into four pieces:
\begin{align*}
S_1=\sum_{i=0}^n\sum_{j=0}^n(\cdot),\quad
S_2=\sum_{i=0}^n\sum_{j=n+1}^{2n}(\cdot),\quad
S_3=\sum_{i=n+1}^{2n}\sum_{j=0}^n(\cdot),\quad
S_4=\sum_{i=n+1}^{2n}\sum_{j=n+1}^{2n}(\cdot).
\end{align*}
For ${2i\choose i}\equiv 0\pmod{p}$ for $n+1\le i\le 2n$, we have
$S_4\equiv 0\pmod{p}$. By the symmetry $i\leftrightarrow j$, we get $S_2=S_3$.
It follows that
\begin{align}
\sum_{i=0}^{p-1}\sum_{j=0}^{p-1}\frac{{2i\choose i}{2j\choose j}}{{i+j\choose i}}\equiv S_1+2S_2\pmod{p}.\label{b-2}
\end{align}

Note that for $0\le i\le n$,
\begin{align}
{2i\choose i}=(-4)^i{-\frac{1}{2}\choose i}\equiv (-4)^i{n\choose i}\pmod{p}.\label{b-3}
\end{align}
Thus,
\begin{align}
S_1&\overset{\eqref{b-3}}{\equiv} \sum_{i=0}^n\sum_{j=0}^n(-4)^{i+j}\frac{ {n\choose i}{n\choose j}}{ {i+j\choose i}}\pmod{p}\notag\\[5pt]
&\overset{\eqref{b-1}}{=}\frac{(-3)^n(2n-1)}{4(n+1)}+\frac{4^n}{{2n\choose n}}
\left(\frac{1}{2}-\sum_{k=0}^nC_k\left(-\frac{3}{4}\right)^{k+1}\right)\notag\\[5pt]
&\equiv -\left(\frac{p}{3}\right)+\frac{(-1)^n}{2}
-(-1)^n\sum_{k=0}^nC_k\left(-\frac{3}{4}\right)^{k+1}\pmod{p},\label{b-4}
\end{align}
where we utilize ${2n\choose n}\equiv (-1)^n\pmod{p}$ in the last step.

Since $C_k\equiv 0\pmod{p}$ for $n+1\le k\le 2n-1$, we have
\begin{align*}
\sum_{k=0}^nC_k\left(-\frac{3}{4}\right)^{k+1}
&\equiv \sum_{k=0}^{2n}C_k\left(-\frac{3}{4}\right)^{k+1}-C_{2n}\left(-\frac{3}{4}\right)^{2n+1}\\ &\overset{\eqref{mt-2}}{\equiv}\frac{1-4^{\frac{p+1}{2}}}{2}-\left(-\frac{3}{4}\right)^p-C_{p-1}\left(-\frac{3}{4}\right)^p
\pmod{p}.
\end{align*}
Using the Fermat's little theorem and
\begin{align*}
C_{p-1}=\frac{{2p-2\choose p-1}}{p}=\frac{{2p-1\choose p-1}}{2p-1}\equiv -1\pmod{p},
\end{align*}
we arrive at
\begin{align}
\sum_{k=0}^nC_k\left(-\frac{3}{4}\right)^{k+1}
\equiv -\frac{3}{2}\pmod{p}.\label{b-5}
\end{align}
Substituting \eqref{b-5} into \eqref{b-4} gives
\begin{align}
S_1&\equiv
2(-1)^n-\left(\frac{p}{3}\right)\pmod{p}.\label{b-6}
\end{align}

Note that
\begin{align}
S_2=\sum_{i=0}^{n}\sum_{j=n+1}^{2n}\frac{{2i\choose i}{2j\choose j}}{{i+j\choose i}}
=\sum_{i=0}^{n}\sum_{j=1}^{n}\frac{{2i\choose i}{2j+2n\choose j+n}}{{i+j+n\choose i}}.\label{b-7}
\end{align}
For $i+j\le n$ and $1\le j\le n$,
\begin{align*}
\frac{{2j+2n\choose j+n}}{{i+j+n\choose i}}\equiv 0\pmod{p},
\end{align*}
and so the summand on the right-hand side of \eqref{b-7} is congruent to $0$ modulo $p$.

On the other hand, for $i+j\ge n+1$ and $1\le j\le n$,
\begin{align}
\frac{{2j+2n\choose j+n}}{{i+j+n\choose i}}&=\frac{i!}{(n+j)!}\cdot\frac{(2n+2)\cdots(2n+2j)}{(2n+2)\cdots(i+j+n)}\notag\\
&\equiv \frac{i!}{(n+j)!}\cdot\frac{(2j-1)!}{(i+j-n-1)!}\pmod{p}.\label{b-8}
\end{align}
It follows from \eqref{b-3} and \eqref{b-8} that
\begin{align*}
\frac{{2i\choose i}{2j+2n\choose j+n}}{{i+j+n\choose i}}\equiv
\frac{(-4)^i{j-1\choose n-i}{2j\choose j}}{2{n+j\choose j}}\pmod{p}.
\end{align*}
Since
\begin{align*}
{n+j\choose j}\equiv {-\frac{1}{2}+j\choose j}=\frac{{2j\choose j}}{4^j}\pmod{p},
\end{align*}
we have
\begin{align}
\frac{{2i\choose i}{2j+2n\choose j+n}}{{i+j+n\choose i}}\equiv
\frac{(-1)^i\cdot4^{i+j}\cdot{j-1\choose n-i}}{2}\pmod{p}.\label{b-9}
\end{align}
Substituting \eqref{b-9} into \eqref{b-7} gives
\begin{align*}
S_2&\equiv \frac{1}{2}\sum_{j=1}^{n}4^j\sum_{i=0}^{n}(-4)^i{j-1\choose n-i}\pmod{p}\\
&=\frac{(-4)^n}{2}\sum_{j=1}^{n}4^j\sum_{i=0}^{n}{j-1\choose i}\left(-\frac{1}{4}\right)^i\quad \text{($i\to n-i$)}\\
&=2(-4)^n\sum_{j=1}^{n}3^{j-1}\\
&=(-12)^n-(-4)^n.
\end{align*}
Thus,
\begin{align}
S_2\equiv \left(\frac{p}{3}\right)-(-1)^n\pmod{p}.\label{b-10}
\end{align}
The proof of \eqref{conj-1} follows from \eqref{b-2}, \eqref{b-6} and \eqref{b-10}.
\qed

\section{Proof of Theorem \ref{t-2}}
\begin{lem}
For any non-negative integer $n$, we have
\begin{align}
\sum_{i=0}^n\sum_{j=0}^n(-4)^{i+j}(i+j)\frac{ {n\choose i}{n\choose j}}{ {i+j\choose i}}=16n(-3)^{n-1}+\frac{8n 4^n}{{2n\choose n}}\sum_{k=0}^n{2k\choose k}\left(-\frac{3}{4}\right)^k.\label{c-1}
\end{align}
\end{lem}
{\noindent \it Proof.}
By the multi-Zeilberger algorithm \cite{az-aam-2006}, we obtain the recurrence for the left-hand side of \eqref{c-1}:
\begin{align*}
&-6(n+1)(581n+793)s(n)+(818n^2-6653n-9936)s(n+1)\\
&+(2166n^2+3474n+2898)s(n+2)+(2n+5)(251n+92)s(n+3)=0.
\end{align*}
It is easy to verify that the right-hand side of \eqref{c-1} also satisfies the above recurrence and both sides of \eqref{c-1} are equal for $n=0,1,2$.
\qed

{\noindent \it Proof of \eqref{a-1}.}
Let $n=\frac{p-1}{2}$. In a similar way,
\begin{align}
\sum_{i=0}^{p-1}\sum_{j=0}^{p-1}(i+j)\frac{{2i\choose i}{2j\choose j}}{{i+j\choose i}}\equiv S_1+2S_2\pmod{p},\label{c-new}
\end{align}
where
\begin{align*}
S_1=\sum_{i=0}^{n}\sum_{j=0}^{n}(i+j)\frac{{2i\choose i}{2j\choose j}}{{i+j\choose i}},
\end{align*}
and
\begin{align*}
S_2=\sum_{i=0}^{n}\sum_{j=n+1}^{2n}(i+j)\frac{{2i\choose i}{2j\choose j}}{{i+j\choose i}}.
\end{align*}

By \eqref{b-3} and \eqref{c-1}, we have
\begin{align}
S_1&\overset{\eqref{b-3}}{\equiv} \sum_{i=0}^n\sum_{j=0}^n(-4)^{i+j}(i+j)\frac{ {n\choose i}{n\choose j}}{ {i+j\choose i}}\pmod{p}\notag\\
&\overset{\eqref{c-1}}{=}16n(-3)^{n-1}+\frac{8n 4^n}{{2n\choose n}}\sum_{k=0}^n{2k\choose k}\left(-\frac{3}{4}\right)^k\notag\\
&\equiv \frac{8}{3}\left(\frac{p}{3}\right)-4(-1)^n\sum_{k=0}^n{2k\choose k}\left(-\frac{3}{4}\right)^k\pmod{p},\label{c-2}
\end{align}
where we make use of ${2n\choose n}\equiv (-1)^n\pmod{p}$ in the last step.

Since ${2k\choose k}\equiv 0\pmod{p}$ for $n+1\le k\le 2n$, we have
\begin{align}
\sum_{k=0}^n{2k\choose k}\left(-\frac{3}{4}\right)^k\equiv \sum_{k=0}^{2n}{2k\choose k}\left(-\frac{3}{4}\right)^k
\overset{\eqref{mt-1}}{\equiv} 4^{\frac{p-1}{2}}\equiv 1\pmod{p}.\label{c-3}
\end{align}
Substituting \eqref{c-3} into \eqref{c-2} gives
\begin{align}
S_1\equiv \frac{8}{3}\left(\frac{p}{3}\right)-4(-1)^n\pmod{p}.\label{c-4}
\end{align}

On the other hand, by \eqref{b-9} we have
\begin{align}
S_2&=\sum_{i=0}^{n}\sum_{j=n+1}^{2n}(i+j)\frac{{2i\choose i}{2j\choose j}}{{i+j\choose i}}\notag\\
&=\sum_{i=0}^{n}\sum_{j=1}^{n}(i+j+n)\frac{{2i\choose i}{2j+2n\choose j+n}}{{i+j+n\choose i}}\notag\\
&\overset{\eqref{b-9}}{\equiv} \frac{1}{2}\sum_{j=1}^{n}4^j\sum_{i=0}^{n}(-4)^i(i+j+n){j-1\choose n-i}\pmod{p}\notag\\
&=\frac{(-4)^n}{2}\sum_{j=1}^{n}4^j\sum_{i=0}^{n}(2n+j-i){j-1\choose i}\left(-\frac{1}{4}\right)^i,\label{c-5}
\end{align}
where we set $i\to n-i$ in the last step.
Note that
\begin{align}
&\sum_{i=0}^{n}(2n+j-i){j-1\choose i}\left(-\frac{1}{4}\right)^i\notag\\
&=(2n+j)\left(\frac{3}{4}\right)^{j-1}-(j-1)\sum_{i=1}^{n}{j-2\choose i-1}\left(-\frac{1}{4}\right)^i\notag\\
&=(2n+j)\left(\frac{3}{4}\right)^{j-1}+\frac{j-1}{4}\left(\frac{3}{4}\right)^{j-2}.\label{c-6}
\end{align}
Substituting \eqref{c-6} into \eqref{c-5} and making elementary calculation gives
\begin{align*}
S_2\equiv \frac{(-12)^n(10n-3)+(-4)^n(3-6n)}{3}\pmod{p}.
\end{align*}
It follows that
\begin{align}
S_2\equiv 2(-1)^n-\frac{8}{3}\left(\frac{p}{3}\right)\pmod{p}.\label{c-7}
\end{align}
Combining \eqref{c-new}, \eqref{c-4} and \eqref{c-7}, we complete the proof of \eqref{a-1}.
\qed


\begin{thebibliography}{99}
\small \setlength{\itemsep}{-.8mm}

\bibitem{ag-jis-2014}E. Allen and I. Gheorghiciuc,
A weighted interpretation for the super Catalan numbers,
J. Integer Seq. 17 (2014), Article 14.10.7, 9 pp.

\bibitem{apagodu-injt-2018}M. Apagodu, Elementary proof of congruences involving sum of binomial coefficients, Int. J. Number Theory, online, doi: 10.1142/S1793042118500938.

\bibitem{az-aam-2006}M. Apagodu and D. Zeilberger, Multi-variable Zeilberger and Almkvist--Zeilberger algorithms and the sharpening of Wilf--Zeilberger theory, Adv. in Appl. Math. 37 (2006), 139--152.

\bibitem{cw-a-2012}X. Chen and J. Wang, The super Catalan numbers $S(m,m+s)$ for $s\le 4$, preprint, 2012, arXiv:1208.4196.

\bibitem{gessel-jsc-1992}I. Gessel, Super ballot numbers, J. Symbolic Comput. 14 (1992) 179--194.

\bibitem{gz-aam-2010}V.J.W. Guo and J. Zeng, Some congruences involving central $q$-binomial coefficients, Adv. in Appl. Math. 45 (2010), 303--316.

\bibitem{mt-jnt-2018}S. Mattarei and R. Tauraso, From generating series to polynomial congruences, J. Number Theory, 182 (2018), 179--205.

\bibitem{ps-rsa-2006}N. Pippenger and K. Schleich, Topological characteristics of random triangulated surfaces, Random Structures Algorithms 28 (2006), 247--288.

\bibitem{stanley-b-1999}R.P. Stanley, Enumerative Combinatorics, Vol. 2, Cambridge Univ. Press, Cambridge, 1999.

\bibitem{st-aam-2010}Z.-W. Sun and R. Tauraso, New congruences for central binomial coefficients, Adv. in Appl. Math. 45 (2010), 125--148.

\bibitem{st-ijnt-2011}Z.-W. Sun and R. Tauraso, On some new congruences for binomial coefficients, Int. J. Number Theory 7 (2011), 645--662.

\bibitem{tauraso-aam-2012}R. Tauraso, $q$-Analogs of some congruences involving Catalan numbers, Adv. in Appl. Math. 48 (2012), 603--614.

\end{thebibliography}
\end{document}